  \documentclass[letterpaper, 10 pt, conference]{IEEEconf}
\IEEEoverridecommandlockouts
\usepackage{cite}
\usepackage{amsmath,amssymb,amsfonts}
\usepackage{graphicx}
\usepackage{textcomp}
\usepackage{xcolor}
\usepackage{soul}
\usepackage{mathtools}
\usepackage{nicefrac}
\usepackage{algorithm}
\usepackage{algpseudocode}
\usepackage{booktabs}
\usepackage{multirow}
\usepackage{setspace} 
\usepackage{bm}
\usepackage{url}

\newtheorem{remark}{Remark}

\algdef{SE}[SUBALG]{Indent}{EndIndent}{}{\algorithmicend\ }%
\algtext*{Indent}
\algtext*{EndIndent}

\def\BibTeX{{\rm B\kern-.05em{\sc i\kern-.025em b}\kern-.08em
    T\kern-.1667em\lower.7ex\hbox{E}\kern-.125emX}}

\usepackage{placeins}
\newcommand{\subparagraph}{}
\usepackage{titlesec}

\begin{document}
\bstctlcite{IEEE:BSTcontrol} % To activate IEEEtran.bst controls in .bib file

\title{
Value-Oriented Forecast Combinations for Unit Commitment
% Evaluating Task-Based Forecast For Power Systems With Wind Energy
% \\
% {\footnotesize }
\thanks{The authors wish to thank the team of the Rutgers Center for Ocean Observation and Leadership (RUCOOL) for their help with data preparation.}
}

\author{Mehrnoush Ghazanfariharandi, Robert Mieth\\
\textit{Industrial and Systems Engineering Department}, 
\textit{Rutgers University},
NJ, USA \\
\{mehrnoush.ghazanfariharandi, robert.mieth\}@rutgers.edu
% \and
% \IEEEauthorblockN{ Robert Mieth}
% \IEEEauthorblockA{\textit{Industrial and Systems Engineering Department} \\
% \textit{Rutgers University}\\
% NJ, USA \\
% Robert.mieth@rutgers.edu}
}

\maketitle

\begin{abstract}
Value-oriented forecasts for two-stage power system operational problems have been demonstrated to reduce cost, but prove to be computationally challenging for large-scale systems because the underlying optimization problem must be internalized into the forecast model training. 
Therefore, existing approaches typically scale poorly in the usable training data or require relaxations of the underlying optimization.
This paper presents a method for value-oriented forecast combinations using progressive hedging, which unlocks high-fidelity, at-scale models and large-scale datasets in training.
We also derive one-shot training model for reference and study how different modifications of the training model impact the solution quality. 

\end{abstract}

% \begin{IEEEkeywords}

% \end{IEEEkeywords}

\section{Introduction}
Operational constraints of large-scale power plants require grid operators to decide on generator schedules ahead of time when load demand and weather-dependent generation are still uncertain.
Hence, these scheduling decisions are made using forecasts of any uncertain quantities.
However, ``good'' forecasts measured by how well they match the true outcome do not necessarily lead to better decisions \cite{carriere2019integrated,stratigakos2024decision}. 
This aspect of forecasting is particularly acute in problems with asymmetric cost functions, as in critical infrastructure like power systems where the cost of resource shortage far exceed the cost of overage \cite{morales2023prescribing}.
This observation has motivated research on \textit{value-oriented} 
forecasting, where the quality of a forecast is measured the decision value\cite{donti2017task,morales2023prescribing,zhang2024toward,dias2025application}.

Value-oriented forecasts offer a practical pathway to internalize more available information into the decision without modifying the decision-making process itself, e.g., through probabilistic forecasts and stochastic programming.
(See also the discussion in \cite{mieth2024prescribed}.)
Essentially, value-oriented forecasts are biased such that they improve the decision without altering the structure of the forecast and the decision-making problem itself.
Value-oriented forecasting models can be trained by unifying the forecast model training problem and the decision-making problem as a bilevel program \cite{zhang2024toward,dias2025application,morales2023prescribing} or by integrating the decision-making problem into a gradient-descent training pipeline \cite{donti2017task}.

Existing approaches as in\cite{donti2017task,dias2025application,zhang2024toward,morales2023prescribing} suffer from poor scalability of the training problem and either remain small-scale or achieve practical scale through intricate heuristics or by giving up modeling details.
In this paper, we obtain value-oriented forecasts at scale even with high-fidelity models in the training phase using progressive hedging (PH). 
In particular, we take the perspective of a power system operator that has access to point forecasts of uncertain demand and renewable generation from multiple forecasting services and seeks a value-oriented combination of these forecasts to achieve lower cost in a two-stage unit commitment problem.

{\color{black}
We propose a scalable value-oriented training of forecast combinations using PH along with a novel modification that we call ``push-forward PH'' and that further improves upon PH computational time.
We use real-world data sets and a power system testbed of practical scale to demonstrate the effectiveness of our method.
}
Relative to similar work in \cite{zhang2024toward,morales2023prescribing} our approach enables the use of more historical data and avoids model approximations {\color{black} and reformulations} during training.
{\color{black}
To compare our approach with standard approaches from the literature, we also formulate and solve the training problem as a single level-equivalent reformulation of a bilevel program.}
Because our approach allows training with high-fidelity models {(\color{black}i.e., a unit commitment formulation with binary variables and network constraints)} we can analyze the impact of relaxing the decision-making problem in training.

\section{Problem description}

We consider a power system operator that manages a set of assets including wind turbines, large-scale generators, and flexible energy resources.
To accommodate the dispatch lead time needed for large-scale generation units, the operator solves a two-stage problem. 
First, the operator solves a unit commitment (UC) problem on the day before the scheduled power delivery to accommodate planning lead times of some generators. 
Then, closer to actual power delivery, the operator solves a second ``real-time'' (RT) problem that uses the previously scheduled units and updated information on demand and renewable injection.

\subsection{Unit commitment and real-time problem}

The power system is modeled as a graph network with a set of nodes $i\in[N]$ (we write $[N]=\{1,..,N\}$) and lines (edges) $l\in[L]$. Every day $d$, the system operator first schedules production $p_{g,t,d}$ and commitment status $u_{g,t,d}$ for each timestep $t\in[T]$ and each generator $g\in[G]$ . 
Schedules depend on uncertain net-load $L_{t,i,d}$ (i.e., load demand minus renewable injection) at time $t$ on day $d$ for node $i$. 
In the day-ahead UC the system operator uses a forecast  $\hat{L}_{t,i,d}$.
We collect all $\hat{L}_{t,i}$ in the matrix $\hat{\bm{L}}$ and write
the resulting 
mixed-integer linear program for day $d$:
\allowdisplaybreaks
\begin{subequations}
\label{eq:mainderetministicuc}
\begin{align}
&  \hspace{-0.93cm} {\rm UC}\big(\hat{\bm{L}}_d\big): \nonumber\\
\min \ 
& \!\!\sum_{g \in [G]} \! \sum_{t=2}^{T} \!\! \left( c_g^{SU} y_{g,t,d} \!+ \!c_g^{SD} (u_{g,t-1,d} -\! u_{g,t,d} +\! y_{g,t,d}) \right)\!+ \nonumber \\ & \!\! \sum_{g \in [G]} \! \sum_{t=1}^{T}\! c_g p_{g,t,d} \!+\!\! \sum_{i \in [N]}\! \sum_{t=1}^{T} \! (c_i^{\rm shed} \hat{l}_{i,t,d} \!+\! c_i^{\rm cur} \hat{w}_{i,t,d} ) \label{objectiveucrt}\\
    & \hspace{-2em}  \text{s.t.}\ \sum_{g\in [G]_i} \!\!p_{g,t,d}  + \!\!\!\!\sum_{l|r(l)=i} \!\!\!\hat{f}_{l,t,d} -\sum_{l|o(l)=i} \hat{f}_{l,t,d} =   \hat{L}_{i,t,d}  - \hat{l}_{i,t,d} \nonumber\\& \quad \ \forall t\in [T] ,\ \forall i \in  [N] \label{eq:UC_balance} \\
      & \hspace{-2em} \hat{f}_{l,t,d} = B_{l} (\hat{\theta}_{o(l),t,d} - \hat{\theta}_{r(l),t,d}) \quad \forall t\! \in\! [T],\ \forall l \!\in \![L] \label{eq:UC_flow1}\\
   &\hspace{-2em} \hat{\theta}_{ref,t,d} = 0 \quad \forall t \in [T]  \label{eq:UC_voltageangle}\\
    &\hspace{-2em} -\overline{F}_{l} \le \hat{f}_{l,t,d} \le \overline{F}_{l} \quad \forall t \in [T] ,\ \forall l \in [L] \label{eq:UC_flow2}\\
& \hspace{-2em} \sum_{i=t-\overline{\ell}_g+1}^{t}\!\!\! y_{g,t,d} \leq u_{g,t,d} \quad \forall t \in [\overline{\ell}_g + 1, T],\ \forall g \in [G] \label{eq:UC_UC1}\\
  & \hspace{-2em} \sum_{i=t-\underline{\ell}_g+1}^{t} \!\!\!\!\!\! y_{g,t,d} \! \!\leq \!1 \!- \!u_{g,t-\underline{\ell}_g,d} \quad \forall t \!\in \! [\underline{\ell}_g \!+ \!1, T],\ \forall g \!\in \! [G] \label{eq:UC_UC2}\\
  & \hspace{-2em} -u_{g,t-1,d} + u_{g,t,d}  \leq y_{g,t,d} \quad \forall t \in [2, T],\ \forall g \!\in \! [G] \label{eq:UC_UC3} \\
  & \hspace{-2em} \underline{P}_g u_{g,t,d} \leq p_{g,t,d} \leq \overline{P}_g u_{g,t,d} \quad \forall t \in [1, T],\ \forall g \in  [G] \label{eq:UC_powerlimit}\\
  & \hspace{-2em}  p_{g,t,d} - p_{g,t-1,d} \!\leq \!R_g u_{g,t-1,d} + \overline{R}_g (1 - u_{g,t-1,d}) \nonumber \\ & \quad \forall t \in [2, T],\ \forall g \in  [G]  \label{eq:UC_ramping1}\\
& \hspace{-2em}    p_{g,t-1,d} - p_{g,t,d} \leq R_g u_{g,t,d} + \overline{R}_g (1 - u_{g,t,d})  \nonumber \\& \quad \forall t \in [2, T],\ \forall g \in  [G] \label{eq:UC_ramping2}\\
% & \hspace{-2em} \hat{w}^{cur}_{i,t} \ge 0\\
& \hspace{-2em} 0 \le \hat{l}_{i,t,d} \le \max\{0,\!\hat{L}_{i,t,d}\!-\!\hat{W}_{i,t,d}\} \label{eq:UC_shed}\  \forall i\!\in\! [N] ,\ \forall t\! \in\! [T]\\
& \hspace{-2em} 0 \le \hat{w}_{i,t,d} \le \hat{W}_{i,t,d} \label{eq:UC_cur} \quad \forall i\in [N] ,\ \forall t \in [T]\\
& \hspace{-2em} u_{g,t,d}  \in\{0,1\}  \qquad \forall g\in [G],\ \forall t \in [2,T]\label{eq:uasbinary}
\end{align}%
\end{subequations}%
\allowdisplaybreaks[0]%
Objective \eqref{objectiveucrt} minimizes the total generator production and start-up/shut-down cost (parametrized by $c_g, c_g^{SU}, c_g^{SD}$) and cost of load shedding $\hat{l}_{i,t,d}$ and curtailment $\hat{w}_{i,t,d}$ (parametrized by $c^{\rm shed}, c^{\rm cur}$). Eq.~\eqref{eq:UC_balance} ensures that scheduled production meets the net demand forecast minus potential load shedding $\hat{l}_{i,t,d}$.
Eq.~\eqref{eq:UC_flow1} models the power flow $\hat{f}_{l,t,d}$ over each line $l$ as a function of line susceptance $B_l$ and the difference between the voltage angle $\hat{\theta}_{i,t,d}$ at the node at the originating end of line $l$, i.e., $o(l)$, and the node at the receiving end of line $l$, i.e., $r(l)$.
Constraint \eqref{eq:UC_voltageangle} defines the voltage angle of the reference node and constraint (\ref{eq:UC_flow2})  transmission capacity limits $\overline{F}_l$.
Constraints~\eqref{eq:UC_UC1}, \eqref{eq:UC_UC2} ensure generator minimum uptime $\overline{\ell}_g$ and downtime $\underline{\ell}_g$. In \eqref{eq:UC_UC3}, $y_{g,t,d}$ captures generator startup.
Constraints \eqref{eq:UC_powerlimit}--\eqref{eq:UC_ramping2} enforce lower and upper limits on power production ($\underline{P}_g$, $\overline{P}_g$) and ramping limits ($R_g$ when online; $\overline{R}_g$ when starting up) for each generator $g$ depending on the binary unit commitment decision $u_{g,t,d}\in\{0,1\}$. 
Lastly, constraints (\ref{eq:UC_shed}) and (\ref{eq:UC_cur}) establish the upper and lower limits on load shedding and renewable curtailment, respectively.
For \eqref{eq:UC_shed} and \eqref{eq:UC_cur} we assume w.l.o.g. that the system operator can separate renewable forecasts $\hat{W}_{i,t,d}$ from the net-load forecasts.

\subsubsection{Real-time problem}
During real-time operations, the uncertain net-load materializes as $\bar{L}_{i,t,d}$ and the system operator manages energy imbalances that result from inaccurate forecasts by solving:
\allowdisplaybreaks
\begin{subequations}
\label{eq:mainderetministicrt}
\begin{align}
&  \hspace{-0.99cm} {\rm RT}\big(\bm p^*_d, \bm u^*_d,  \bar{\bm L}_d\big): \nonumber\\
\min \ 
&\sum_{g \in [G]} \sum_{t=1}^{T} (c_g^+ r_{g,t,d}^+ + c_g^- r_{g,t,d}^- )+ \sum_{i \in [N]} \sum_{t=1}^{T} (c_i^{\rm shed} l_{i,t,d} \nonumber \\ &  +  c_i^{\rm cur} w_{i,t,d})  \label{objectivert}\\
    & \hspace{-2em}  \text{s.t.}\  \sum_{g\in [G]_i} \!\!( p^*_{g,t,d}+  r_{g,t,d}^{+} - r_{g,t,d}^{-}) + \!\!\!\sum_{l|r(l)=i} \!\!\!f_{l,t,d} -\!\!\!\sum_{l|o(l)=i} \!\!f_{l,t,d} \nonumber\\& =   \bar{L}_{i,t,d} - l_{i,t,d} \quad \forall t\in [T],\ \forall i \in  [N] \label{eq:RT_balance}\\
    &\hspace{-2em} f_{l,t,d} = B_{l} (\theta_{o(l)} - \theta_{r(l)}) \quad \forall t \in [T] ,\ \forall l \!\in \![L] \label{eq:RT_flow1} \\
   &\hspace{-2em}\theta_{ref,t,d} = 0 \quad \forall t \in [T] \label{eq:RT_voltageangle}\\
    & \hspace{-2em}-\overline{F}_{l} \le f_{l,t,d} \le \overline{F}_{l} \quad \forall t \in [T] ,\ \forall l \in [L] \label{eq:RT_flow2}\\
  &\hspace{-2em} \underline{P}_g u^*_{g,t,d} \!\!\leq p^*_{g,t,d} \!+ \!r^+_{g,t,d} \!- r^-_{g,t,d} \!\! \leq \!\! \overline{P}_g u^*_{g,t,d} \quad \forall t\!\! \in [1, T], \nonumber\\& \forall g\in [G] \label{eq:RT_powerlimit}\\
  & \hspace{-2em}  (p^*_{g,t,d} \!+\! r^+_{g,t,d}\!\!- r^-_{g,t,d})\! - \!(p^*_{g,t-1,d} \!+ r^+_{g,t-1,d} -\! r^-_{g,t-1,d} )\nonumber\\& \hspace{-2em}\leq\! R_g u^*_{g,t-1,d} \!+\! \overline{R}_g (\!1 \!-\! u^*_{g,t-1,d}) \! \quad \forall t \!\in \! [2, T],\ \forall g \!\in \!\! [G]\label{eq:RT_ramping1}\\
& \hspace{-2em}  (p^*_{g,t-1,d} \!+ r^+_{g,t-1,d} - r^-_{g,t-1,d})\!\! - (p^*_{g,t,d} + r^+_{g,t,d} - r^-_{g,t,d}) \nonumber \\&\hspace{-2em} \leq \! R_g u^*_{g,t,d} \!+ \!\overline{R}_g (1 - u^*_{g,t,d}) \quad  \forall t \in [2, T],\ \forall g \!\in \! [G] \label{eq:RT_ramping2}\\
    % &\hspace{-2em}  0 \le r_{g,t,d}^+ \le R_g \quad \forall g\in [G],\ \forall t \in [T] \label{eq:RT_up} \\
    & \hspace{-2em} 0 \le r_{g,t,d}^-, r_{g,t,d}^+\le R_g \quad \forall g\in [G],\ \forall t \in [T] \label{eq:RT_down}\\
    &\hspace{-2em} 0 \le l_{i,t,d} \le \max\{0,\!\bar{L}_{i,t,d}\!\!-\!\!\bar{W}_{i,t,d}\} \label{eq:RT_shed}\  \forall i\!\in\! [N] ,\ \forall t\! \in\! [T]\\
    &\hspace{-2em} 0 \le w_{i,t,d} \le \Bar{W}_{i,t,d} \label{eq:RT_cur} \quad \forall i\in [N] ,\ \forall t \in [T]
\end{align}
\end{subequations}
Values $p^*_{g,t,d}$, $u^*_{g,t,d}$ are the decisions obtained from UC. 
For given UC decisions $p^*_{g,t,d}$, $u^*_{g,t,d}$ and net-load realizations $\bar{L}_{i,t,d}$, the RT problem minimizes the cost of upward and downward redispatch ($r^+_{g,t,d}$ and $r^-_{g,t,d}$ with respective cost $c_g^+$ and $c_g^-$) such that the power balance in Eq.~\eqref{eq:RT_balance} is ensured.
Constraints \eqref{eq:RT_powerlimit}--\eqref{eq:RT_ramping2} limit redispatch based on $p^*_{g,t,d}$, $u^*_{g,t,d}$ and the generation and ramping limits. 
The remaining constraints \eqref{eq:RT_flow1}--\eqref{eq:RT_cur} are functionally equivalent to their analogous constraints in \eqref{eq:mainderetministicuc}. 

\subsection{Optimal forecast combination}

To solve UC, the system operator requires a single net-load forecast.
We assume that the operator has access to  net-load forecasts from multiple providers and must decide how to combine the information from these forecasts.
We denote the combined net-load forecast as $\hat{L}_{i,t,d}^{\rm comb} =\!\sum_{k=1}^{K} \!\lambda_{k} \hat{L}_{i,t,d,k}$, where $K$ is the number of forecast providers and $\lambda_{k}$ is a provider-specific weight.

Typically, forecast quality is measured by how well it matches the observed value and numerous statistical methods exist to this end \cite{wang2023forecast}. 
However, from the perspective of the two-stage UC and RT problem, a better metric for a good forecast combination, i.e., the choice of $\lambda_k$, is related to the cost of system operation resulting from the daily forecast-observation pair.
Formally:
\allowdisplaybreaks
\begin{subequations}
\begin{align}
\bm{\lambda}^* \!\!= \arg\min_{\bm{\lambda}} \ 
    & \mathbb{E}_{\bm{L}}\!\!\left[\!{\rm UC}\!\left( \sum_{k=1}^{K} \lambda_k \hat{\bm{L}}_{k} \!\!\right) \!+ {\rm RT}\left(\bm{p}^*, \bm{u}^*, \bm{L} \right)\right] \\
\text{s.t.} \quad 
    & \bm{p}^*, \bm{u}^* \in \arg\min UC\left( \sum_{k=1}^{K} \lambda_k \hat{\bm{L}}_{k} \right).
\end{align}%
\label{eq:forecast_comb_general}%
\end{subequations}%
\allowdisplaybreaks[0]%
Bold variables indicate vectors, which we use to omit some indices for easier readability.
\label{sssec:training_problem}
To determine a forecast combination that achieves the desired property in \eqref{eq:forecast_comb_general}, we use historical forecasts and actual observations and formulate the bilevel training program:
\allowdisplaybreaks
\begin{subequations}
\label{eq:mainbilevel}
 \begin{align}
   \hspace{0.02em} \min_{\lambda_1,...,\lambda_K} \
&  \frac{1}{D}(\sum_{d=1}^D \{\eqref{objectiveucrt} + \eqref{objectivert}\}) \label{ew:bilevel_obj}\\
& \hspace{-2.6em}  \text{s.t.}\ \hat{\bm{L}}^{comb}_d\! =\!\! \sum_{k=1}^{K} \!\lambda_{k} \hat{\bm{L}}_{k,d} \ , \ {\sum_{k=1}^{K} \lambda_{k}} =1 \label{eq: forecast combination}\\
& \hspace{-1em}\eqref{eq:RT_balance}\!-\!\eqref{eq:RT_cur}\ [\text{Real-time constraints}] \ \forall d\! \in \! [D]\label{eq:secondstageinbilevel}\\
& \hspace{-0.9em} \bm{p}^*_{d}, \bm{u}^*_{d} \in \!\! \Bigg\{ \text{arg}\min_{\bm{p}_{d}, \bm{u}_{d}}  \eqref{objectiveucrt} 
\nonumber\\
    & \hspace{0em}  \text{s.t.}\!\!\!\ \sum_{g\in [G]_i}\!\! p_{g,t,d}  +\!\! \sum_{l|r(l)=i}\! \hat{f}_{l,t,d} -\!\!\sum_{l|o(l)=i} \! \hat{f}_{l,t,d}= \hat{L}^{\rm comb}_{i,t,d}    \nonumber\\& \hspace{1.4em} -\! \hat{l}_{i,t,d} \ \ \forall t\in [T] ,\ \forall i \in  [N] \label{eq: innerprob}\\
 &\hspace{1.4em}\eqref{eq:UC_flow1}\!-\!\eqref{eq:uasbinary} [\text{Day-ahead UC constraints}] \!\!\Bigg\} \forall d\!\! \in\! \![D]. \nonumber 
 \end{align}%
 \end{subequations}%
\allowdisplaybreaks[0]%
Here, the upper-level problem finds a $\bm{\lambda}$ that minimizes the sample average two-stage operation cost over $D$ days for which historical data is available.
Constraint \eqref{eq: forecast combination} computes the combined forecast by assigning weight to each forecast vector and enforces a convex combination of forecasts \cite{stratigakos2024decision}. 
The outer level problem also obtains the RT solution \eqref{eq:secondstageinbilevel} for a given realized net load $\bar{\bm{L}}$ and a given previous unit commitment decision $\bm{p}^*_{d}, \bm{u}^*_{d}$.
The lower-level problem in \eqref{eq: innerprob} solves the UC problem \eqref{eq:mainderetministicuc} and is  parameterized by the combined forecast $\hat{\bm{L}}^{com}_d$.

 {\color{black}
\begin{remark}
There has been extensive research on solving the UC and RT two-stage problem with stochastic methods  to reduce energy cost and increase reliability \cite{roald2023power}.
However, grid operators are unlikely to implement such methods in the near future \cite{mieth2024prescribed}.
Hence, obtaining a point-forecast that improves the decisions made in the \textit{deterministic} UC and RT process, i.e., without changing the UC or RT problem itself, is a central motivation of this paper.
\end{remark}

\begin{remark}
Forecast combination weights can also be generalized to be node- or time-specific, e.g., $\hat{L}^{\text{comb}}_{i,t,d} = \sum_{k=1}^{K} \lambda_{k,i} \hat{L}_{i,t,d,k}$ for node-specific weights. 
However, we only consider provider-specific weights as defined in (3) in this paper to reflect the economic relevance of each provider, enabling the system operator to assess providers based on the usefulness of their forecasts.
\end{remark}
}

\section{Solution Methodology}
Problem \eqref{eq:mainbilevel} is hard to solve not only because it is a bilevel program but also because of the scale of practical power systems, binary variables in the lower-level UC problem, and the scaling with size the  training dataset.
In the following,
{\color{black}
we first discuss an approach to solve \eqref{eq:mainbilevel}  that is common in the literature based on reformulation and then propose an alternative more tractable approach based on decoupling}.

\subsection{{\color{black}Standard approach:} Single-level reformulation}
{\color{black}
A standard approach to solving problem \eqref{eq:mainbilevel} is formulating a single-level equivalent where the lower-level (inner) problem is replaced with its Karush–Kuhn–Tucker (KKT) optimality conditions \cite{zhang2024toward,morales2023prescribing}. 
However, the lower-level UC problem is non-convex due to its binary variables, which obstructs the direct application of the KKT conditions. 
To overcome this, a convex hull relaxation of UC, e.g., as proposed in \cite{hua2016convex}, can be applied. 
We show it as ``UC-R'' in the Appendix below. 
This relaxation allows substituting the KKT conditions of UC-R for the lower level problem in \eqref{eq:mainbilevel}.
The resulting nonlinearities in the complementarity slackness conditions can then be addressed using a regularization approach \cite{scholtes2001convergence} or the traditional Fortuny–Amat ``Big-M'' method. We use the shorthand ``ST-N'' for the resulting single-level equivalent of \eqref{eq:mainbilevel} using the regularization approach, and ``ST-M'' for its counterpart using the Big-M method.
We omit these standard formulations here, but, because we use them for comparison to our proposed method below, we report all detailed formulations for reproducibility of ST-N and ST-M as supplementary material at \cite{ghazanfariharandi2025value}.}

The single-level equivalent problems solve a {\color{black}large-scale} two-stage network-constrained problem for each time $t$ over all days $d$ and include binary (ST-M) or non-linear (ST-N) structures.
Therefore, they cannot be expected to be computationally tractable for practical application. 
Next, we propose a more tractable solution approach using decomposition.
{\color{black}
\begin{remark}
Relaxing the inner UC problem with its convex hull provides a lower bound for the UC problem. The bilevel program, on the other hand, is generally neither relaxed nor restricted by a relaxed inner problem \cite{kleinert2021survey}.
We therefore note that we describe ST-M and ST-N not as a benchmark for global optimality but for comparison of our proposed approach with standard approaches from the literature. 
\end{remark}
}

 \subsection{{\color{black}Proposed approach:} Progressive hedging algorithm}
 \label{sec:PHA}
Progressive hedging (PH) \cite{rockafellar1991scenarios} decomposes the original problem such that each scenario (day of training data in our case) can be solved independently and then uses an augmented Lagrangian approach to achieve consensus between shared variables. 
This structure makes the algorithm particularly suited for parallelization and drastically reduces the size of each individual problem, allowing for efficient simultaneous computations.
Moreover, because the PH problem solves the two-stage UC and RT problem individually for each day, we do not require the KKT conditions of the UC problem and can instead formulate a combined UC and RT problem using their primal formulations \eqref{eq:mainderetministicuc} and \eqref{eq:mainderetministicrt}:
\allowdisplaybreaks
\begin{subequations}
\label{mainPH}
\begin{align}
&  \hspace{-1.18cm} {\rm PH}\big(\hat{\bm{L}}_d, \bm{\mu}, \rho, \bar{\bm{\lambda}} \big): \nonumber\\
\min_{\bm{\lambda}_d} \
    & \eqref{objectiveucrt} + \eqref{objectivert}+ (\bm{\mu}_d^{\tau-1})^T \bm{\lambda}_d + \frac{\rho}{2} \|\bm{\lambda}_d - \bar{\bm{\lambda}}^{\tau-1}\|^2\label{eq:objtwosatgePH}\\
\text{s.t.} \quad 
    & \hat{\bm{L}}_d^{comb} = {\sum_{k=1}^{K} \lambda_{d,k} \hat{\bm{L}}_{d,k} } \ , \  \sum_{k=1}^K \lambda_{d,k} = 1 \\
        & \sum_{g\in [G]_i} p_{g,t,d}  + \sum_{l|r(l)=i} \hat{f}_{l,t,d} -\sum_{l|o(l)=i} \hat{f}_{l,t,d} \nonumber\\& \hspace{2em}=   \hat{L}^{\rm comb}_{i,t,d}  - \hat{l}_{i,t,d} \quad \ \forall t\in [T] ,\ \forall i \in  [N] \label{eq: bilevellowerlevelbalance}\\
 &\eqref{eq:UC_flow1}-\eqref{eq:uasbinary} \quad[\text{Day-ahead UC constraints}] \label{eq:ph-uc}\\
     & \eqref{eq:RT_balance}-\eqref{eq:RT_cur} \quad[\text{Real-time constraints}]
\end{align}%
\end{subequations}%
\allowdisplaybreaks[0]%
In this formulation, we again use bold symbols to denote vectors. 
{\color{black}We note that the formulation in \eqref{mainPH} formulates exact feasibility of the RT and UC problems but relaxes the lower level optimality condition of the original bilevel problem \eqref{eq:mainbilevel}.
This ensures feasibility for RT and UC, and produces a lower bound on \eqref{eq:mainbilevel} \cite{kleinert2021survey}.
}%r2c2
Problem \eqref{mainPH} is solved for each day in each iteration of the PH algorithm.
In essence, each day computes its individual optimal forecast combination $\bm{\lambda}_d$ based on the data for that day. 
The PH algorithm then computes the average forecast combination $\bar{\bm{\lambda}}$ and each day recomputes its optimal forecast combination with additional PH {\color{black} multiplier} $\bm{\mu}_d$ and {\color{black} penalty} $\rho$ that are added to the objective \eqref{eq:objtwosatgePH}. We can solve \eqref{mainPH} using both UC and UC-R models. To solve \eqref{mainPH} with UC-R, we replace equations \eqref{eq:UC_flow1}--\eqref{eq:uasbinary} in \eqref{eq:ph-uc} with \eqref{eq:UC_flow1}--\eqref{eq:UC_cur}, \eqref{eq:ramprelax1}--\eqref{eq:endrelax} in \eqref{eq:ph-uc}. 
{\color{black}We highlight that RT and UC are feasible for any net-load and net-load forecast, and thus any $\bm \lambda$ produced by the PH algorithm, because of the option to shed load and curtail renewable injection. 
}%R10C5

Alg.~\ref{alg:phmethod} shows the PH method in detail.
After initialization, PH multiplier $\bm{\mu}_d$ is calculated for each training day $d$ based on
the difference between the individual $\bm{\lambda}_d$ and the average $\bar{\bm{\lambda}}$ over all days (Line~\ref{algl:ini_mu} in Alg.~\ref{alg:phmethod}). 
Each day then re-solves \eqref{mainPH} parametrized by the current PH multiplier $\bm{\mu}_d$ and the average $\bar{\bm{\lambda}}$ from the previous solution. (Line~\ref{alogl:solve_prob} in Alg.~\ref{alg:phmethod}). 
These steps are repeated until the consensus gap $g$ (Line~\ref{algl:convergence} in Alg.~\ref{alg:phmethod}) is smaller
than a predefined threshold $\epsilon$.

{\color{black}
\begin{remark}
Since UC is a mixed-integer problem, the PH algorithm cannot guarantee global optimality, but it has been shown to yield high-quality solutions\cite{gade2016obtaining}.
\end{remark}}

\begin{algorithm}
\small
\caption{PH Algorithm}
\label{alg:phmethod}
\begin{algorithmic}[1]
\setstretch{0.85}
 \State \textbf{Input:} $\{\rho >0, \epsilon >0,\hat{\bm{L}}, \bar{\bm{L}}\}$ 
\State \textbf{Initialization:} \label{algl:initstart} 
\Indent
\State $\bm{\lambda}_d^0 \gets  {\rm PH}(\hat{\bm{L}}_d, 0, 0, 0),\ \forall d\in [D]$ \Comment{Solves \eqref{mainPH}}
        \State $\bar{\bm{\lambda}}^0 \gets \frac{1}{D} \sum_{d \in [D]} \bm{\lambda}_d^0$ \label{algl:ini_exp}\Comment{Average weights} 
    \State $\bm{\mu}_d^0 \gets \rho (\bm{\lambda}_d^0 - \bar{\bm{\lambda}}^0),\ \forall d\in [D]$ \label{algl:ini_mu}\Comment{Initial PH multipliers}
\EndIndent \label{algl:initend}
\State $\tau = 1$ \label{algl:initend}\Comment{Set iteration counter}
\Repeat
    \State 
        $\bm{\lambda}_d^\tau \gets {\rm PH}(\hat{\bm{L}}_d, \bm{\mu}_d^{\tau-1}, \rho, \bar{\bm{\lambda}}^{\tau-1}),\ \forall d\in [D]$ \Comment{Solves \eqref{mainPH}} \label{alogl:solve_prob}
    \State $\bar{\bm{\lambda}}^{\tau} \gets \frac{1}{D} \sum_{d \in [D]} \bm{\lambda}_d^\tau$ \label{algl:exp}\Comment{Average weights} \label{algl:average}
    \State $\bm{\mu}_d^\tau \gets \bm{\mu}_d^{\tau-1} + \rho (\bm{\lambda}_d^\tau - \bar{\bm{\lambda}}^\tau),\ \forall d\in [D]$\label{algl:mu} \Comment{\parbox[t]{.2\linewidth}{Current PH\\ multipliers}}
    \State $g^\tau \gets \sum_{d \in [D]} \|\bm{\lambda}_d^\tau - \bar{\bm{\lambda}}^\tau\|$\label{algl:convergence} \Comment{Current convergence}
    \State $\tau \gets \tau + 1$ \label{algl:step}\Comment{Step iteration counter}
\Until{$g^\tau < \epsilon$}
\State {\color{black} \textbf{Output:} $\{ \bm{\bar\lambda^{\tau}}\}$}
\end{algorithmic}
\end{algorithm}
\begin{algorithm}
\small
\caption{Push-forward PH Algorithm (PFPH)}
\label{alg:IMphmethod}
\begin{algorithmic}[1]
\setstretch{0.85}
\State \textbf{Input:} $\{\rho >0, \epsilon >0,\hat{\bm{L}}, \bar{\bm{L}}, D'\}$ 
\State [Lines \ref{algl:initstart}--\ref{algl:initend} of Alg.~\ref{alg:phmethod}]
\Repeat
\State $ds_d \gets \|\bm{\lambda}_d^{\tau-1} - \bar{\bm{\lambda}}^{\tau-1}\|,\ \forall d \in [D]$ \Comment{Deviation score}\label{alog2:solve_ds}
    \State 
    Find $\mathcal{S}', \bar{\mathcal{S}}'$ such that $\mathcal{S'} \subseteq [D],\ |\mathcal{S}'|=D',\ \bar{\mathcal{S}}'=[D]\setminus\mathcal{S}'$ where $\ \forall i\in\mathcal{S}',\ \forall j\in \bar{\mathcal{S}}':\ ds_i \ge ds_j$ \label{alog2:solve_s'}
    \State
        $\bm{\lambda}_d^\tau \gets {\rm PH}(\hat{\bm{L}}_d, \bm{\mu}_d^{\tau-1}, \rho, \bar{\bm{\lambda}}^{\tau-1}),\ \forall d\in \mathcal{S}'$ \Comment{Solves \eqref{mainPH}} \label{alog2:solve_prob}
    \State
        $\bm{\lambda}_d^\tau \gets \bm{\lambda}_d^{\tau-1},\ \forall d\in\bar{\mathcal{S}}'$
    \State [Lines \ref{algl:average}--\ref{algl:step} of Alg.~\ref{alg:phmethod}]
\Until{$g^\tau < \epsilon$}
\State {\color{black} \textbf{Output:} $\{ \bm{\bar\lambda^{\tau}}\}$}
\end{algorithmic}
\end{algorithm}

\subsection{Push-forward progressive hedging}
In each iteration, the PH algorithm re-solves all training days. (See Line~\ref{alogl:solve_prob} in Alg.~\ref{alg:phmethod}.)
{\color{black} We propose modification of Alg.~\ref{alg:phmethod} that speeds up computation by creating a smaller subset $\mathcal{S}' \subseteq [D]$ of the days to be re-evaluated at each iteration (Line~\ref{alog2:solve_prob} in Alg.~\ref{alg:IMphmethod}). This set is obtained as follows:} 
First, we compute a deviation score $ds_d = \|\bm{\lambda}_d^{\tau-1} - \bar{\bm{\lambda}}^{\tau-1}\|$ at each iteration.
{\color{black}
This score reflects how much each day deviates from the consensus. We then select the indices of the $D'$ largest deviation scores to form the subset $S'$. 
}
We also define $\bar{\mathcal{S}}' \coloneqq [D]\setminus \mathcal{S}'$. (See line~\ref{alog2:solve_s'} in Alg.~\ref{alg:IMphmethod}.)
At each iteration, only the $\bm{\lambda}_d$ for $d\in\mathcal{S}'$ are re-computed.
This allows to control of the per-iteration computational cost to achieve fast iterations at a potential longer convergence time. 
{\color{black}Once the total consensus gap g (Line 9 in Alg.~\ref{alg:IMphmethod}) is smaller than
a predefined threshold $\epsilon$, then the algorithm outputs the weights $\bar{\lambda}_k ,\ \forall k \in K$ assigned to each forecast provider.}

\section{Numerical Experiments}

{\color{black}
We test our method using real-world offshore wind data from two different numerical weather prediction models.}

\subsection{Description of experiments and data}\label{ssec:data_description}
We first focus on the IEEE 24-bus test system.
{\color{black}We extend the single-period load data from the original IEEE 24-bus test system to one year of hourly data by scaling it using yearly} real-world demand profiles from ENTSO-E \cite{OpenPowerSystemDataplatform}. For uncertain wind power injections, {\color{black}we use two real-world data sources that include two
forecasts from two different numerical weather simulation models (RUWRF \cite{dicopoulos2021weather,RUCOOL2019} and NOW23 \cite{now23}). }%r11c6
 Each dataset contains hourly day-ahead forecasts and actual measurements
for wind speed, which we translated to wind power using the
NREL 15-Megawatt Reference Wind Turbine \cite{gaertner2020iea}. We locate wind farms at nodes 3, 5, 9, 16,
19, 20 and the capacity of each wind farm is 400 MW. {\color{black}At buses with wind turbines, we define the net load $\hat{\bm L}$ as the load minus the forecasted wind power; At other buses, $\hat{\bm L}$ equals the load demand.} We note we model demand as deterministic and only wind as
uncertain. This, however, is no advantage or disadvantage for the performed computations.
{\color{black} Based on the available data we can model two forecast providers: RUWRF ($\hat{\bm L}_{k=1}$) and NOW23 ($\hat{\bm L}_{k=2}$).} 
We set $c_i^{\rm cur}=\$50/MWh,\ \forall i$ and $c_i^{\rm shed} = \$25,000/MWh,\ \forall i$. 
The time horizon is 24 hours.
We set $D'= \lceil\nicefrac{1}{3}D\rceil$.
{\color{black}
We tune $\rho$ by running Alg.~1 for various values $\rho \in $ [500,\ 1000,\ 5000,\ 10000,\ 15000,\ 20000,\ 25000,\ 30000,\ 35000] and selecting $\rho=25000$ as the one with the lowest testing cost.
}
We set $\epsilon = 10^{-5}$ for PH and PFPH.

All computations have been implemented in Julia using JuMP \cite{lubin2023jump} and solved using Gurobi solver on the Rutgers Amarel cluster on nodes with 128 GB of memory and 16 cores (Dual Intel Xeon Gold 6448Y processors).

\subsection{Small-scale reference cases}
We first solve {\color{black} the single-level training problem} in the  ST-M and ST-N variants to obtain reference solutions and gauge the scalability of this direct approach. 
The largest instance that could be solved without running out of memory or hitting a limit of three computation days used 30 days of historical data. 
We therefore used an instance of the problem with one month worth of training days to compare ST-M and ST-N with the PH methods.
Table~\ref{tab:results} summarizes the results. Values for PH and PFPH have been obtained by solving 
Alg.~\ref{alg:phmethod} and Alg.~\ref{alg:IMphmethod}, respectively, with the UC-R variant of \eqref{mainPH}.
We do this for better comparability because ST-M and ST-N inherently require UC-R.

We observe that the resulting forecast combinations $\lambda_1$ and $\lambda_2$ are similar across the methods and confirm the correct convergence of PH and PFPH. See also the top plot in Fig.~\ref{fig:Lambda_Convergence_Combined_Plot_week}.
Table~\ref{tab:results} also shows the required training time.  
Clearly, the PH and PFPH approaches outperform the ST-M and ST-N methods in terms of computational speed. 
ST-N and PH are similar, but ST-M failed to scale to larger problem instances. 

We test the performance of the obtained forecast combination by running the two-stage UC+RT problem (UC in its standard form with binary variables) for one year with different forecast combinations:
(a) $\lambda_1=1$, $\lambda_2=0$ (using only forecasts from forecast provider 1), (b) $\lambda_1=0$, $\lambda_2=1$ (using only forecasts from forecast provider 2), and (c) $\lambda_1=0.5$, $\lambda_2=0.5$ (using the naive average).
We define the resulting average two-stage cost as ${\rm TST}(\lambda_1, \lambda_2)$.
We denote the average testing results using the value-oriented forecast as ${\rm TST}^*$.
In Table~\ref{tab:results}, columns $\Delta_a$, $\Delta_b$, $\Delta_c$ then show the average daily improvement the value-oriented forecast achieves over the reference methods (a), (b), (c). 
In all cases, ${\rm TST}^*$ improves the solution, indicating systematic benefits of using the value-oriented forecast from any method.

\begin{table}
\centering
 \setlength{\tabcolsep}{2.7pt}
\caption{Results of 4 different solution methods for One Month.} 
\footnotesize
\renewcommand{\arraystretch}{1.2} 
\begin{tabular}{lccccccc}
\hline
\text{Method} & $\lambda_1^*$ &  \text{$\lambda_2^*$}  & \text{Time (s)} & \text{TST* [\$] } & $\Delta_a$ & $\Delta_b$ & $\Delta_c$  \\
\hline
ST-M  & 0.464 & 0.535 &  64800 & 3645874 & -66625 & -58291 & -3171 \\
ST-N & 0.462 & 0.537 & 2368 & 3644656 & -67843 & -59508 & -4388 \\
PH (CR) & 0.471 & 0.528  & 2149 & 3640312 & -72187 & -63853 & -8733 \\
PFPH (CR)& 0.466 & 0.533  & 1943 & 3643742 & -68757 & -60422 & -5302 \\
\hline
\end{tabular}
\begin{flushleft}
\footnotesize
$\Delta_a$: $\!\rm TST^* \!\!\!-\! \rm TST(\!1,\!0\!)$, $\Delta_b$: $\!\rm TST^*\! \!-\! \rm TST(\!0,\!1\!)$, $\Delta_c$: $\rm \!TST^* \!\!\!-\! \rm TST(\!0.5,\!0.5\!)$
\end{flushleft}
\label{tab:results}
\end{table}

\subsection{PH algorithm with full training dataset}

We now investigate modifications of the PH approach using all the available training data (one year). For this dataset, the direct training approaches ST-M and ST-N were intractable and are no longer considered.
We study the following modifications: As above, (PF)PH (CR) solves Algs.~\ref{alg:phmethod}, \ref{alg:IMphmethod} using the UC-R version of \eqref{mainPH}. (PF)PH (B) uses the standard formulation of UC with binary variables as written in \eqref{mainPH} for Algs.~\ref{alg:phmethod}, \ref{alg:IMphmethod}. 
For comparison with \cite{morales2023prescribing}, we also solve a training version without network constraints in the UC stage. This is denoted N (with network) and NN (no network) in Table~\ref{tab:results_binary}, which summarizes the results.

\begin{figure}
    \centering
    \includegraphics[width=0.8\linewidth]{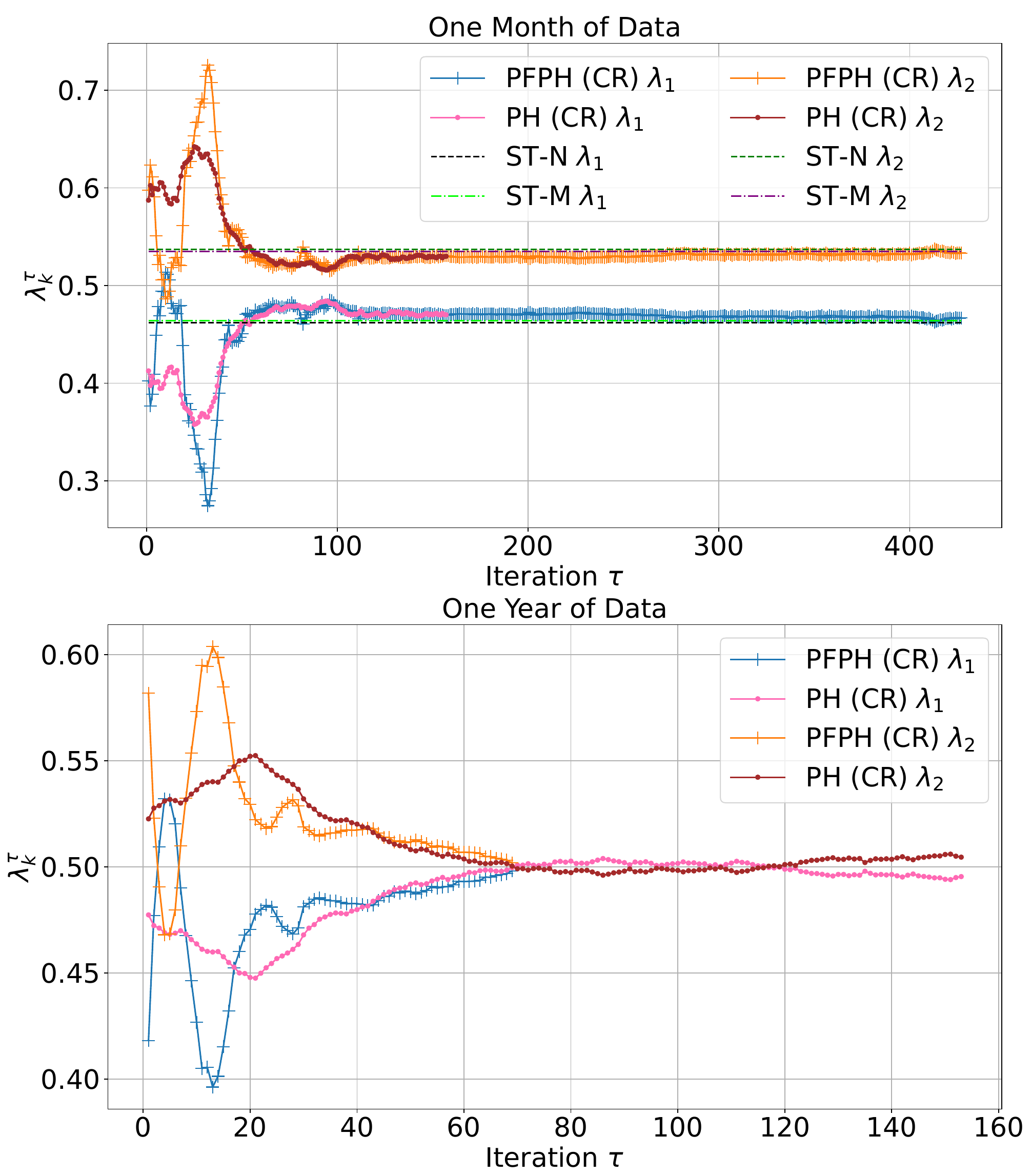}
  \caption{Convergence of $\lambda_k^\tau$ for PH (CR) and PFPH (CR) using one month (top) and one year (bottom) of training data.
  Top plot shows comparison with $\lambda_k$-s obtained from ST-M and ST-N. }\label{fig:Lambda_Convergence_Combined_Plot_week}
\end{figure}

PFPH reduces the solution time of PH without much compromise in the solution, making large-scale problems tractable.
Additionally, Fig.~\ref{fig:Lambda_Convergence_Combined_Plot_week} shows the convergence speed of PH and PFPH in different time periods. 
While PFPH is slightly less stable in earlier iterations, it shows similar convergence behavior as PH. In fact, it meets the convergence criterion faster than PH.
Fewer iterations {\color{black}until convergence} also explain the faster training of PFPH (B) over PFPH (CR), {\color{black} 66 iterations compared to 71 iterations, respectively.}

We observe that forecast combinations strictly improve the decision value. Interestingly, using the convex hull version of the UC problem in training leads to the best improvements in (PF)PH (CR). 
We suspect that the PH algorithm benefits from the convexity of the underlying problem 
{\color{black}as the PH algorithm guarantees a globally optimal solution for convex problems.
(Hence $\bm \lambda$ computed by PH (CR) method is globally optimal in the relaxed space and serves as a lower bound for \eqref{mainPH}.
}Training the forecast combination network-ignorant (rows NN in Table~\ref{tab:results_binary}) improves upon using a single forecast, but performs worse in testing than just averaging the forecasts. 
Here, ignoring network congestion in training creates an advantage for the forecasts from forecast provider 2, which tends to underestimate hourly wind power fluctuations.

\begin{table}
\centering
  \setlength{\tabcolsep}{1.5pt}
\vspace{-6pt}
\caption{Results for one year of data}
\vspace{-6pt}
\footnotesize
\renewcommand{\arraystretch}{1.2}
\begin{tabular}{clccccccc}
\hline
& \text{Method} & $\lambda_1^*$ &  \text{$\lambda_2^*$}  & \text{Time (s)} & \text{TST* [\$] } & $\Delta_a$ & $\Delta_b$ & $\Delta_c$ \\
\hline
\multirow{4}{*}{\rotatebox{0}{N}} 
& PH (B)  & 0.502 & 0.497 &  26357 & 1806985 & -40048 & -45146 & -7671 \\
& PFPH (B) & 0.501 & 0.498 &  3967 & 1805905 & -41128 & -46226 & -8751 \\
& PH (CR) & 0.495 & 0.504 &  11301 & 1805192 & -41841 & -46939 & -9464 \\
& PFPH (CR)& 0.498 & 0.501 &  4225 & 1805590 & -41443 & -46541 & -9066 \\
\hline
\multirow{2}{*}{\rotatebox{0}{NN}} 
& PH (B) & 0.453 & 0.547 & 10956 & 1816795 & -30238 & -35336 & 2139 \\
& PFPH (B)& 0.488 & 0.512 &  4954 &1820322 & -26711 & -31809 & 5666 \\
\hline
\multirow{1}{*}{\rotatebox{0}{}} 
& RMSE & 0.483 & 0.517 & 0 &  1840149& -6884 & -11982 & 25493 \\
\hline
\end{tabular}
\label{tab:results_binary}
\end{table}

\subsection{Comparison to statistical method}
We also compare our value-oriented forecast combination with an established statistics-based combination  method based on the root mean square errors (RMSE) of the forecasts \cite{nowotarski2014empirical}.
For each $k\in[K]$ we compute the RMSE as
\begin{align}
   & \hspace{-0.2cm} \text{RMSE}_k \!=\!\! \left(\!\!\frac{1}{D}\! \sum_{d=1}^{D} (\!\frac{1}{N \times T} \!\sum_{i=1}^{N} \sum_{t=1}^{T} (\hat{L}_{i,t,d,k} -\! \Bar{L}_{i,t,d})^2)\!\!\right)^{\!\frac{1}{2}}\label{eq:rmse}
\end{align}
and then calculate $\lambda_k$ inversely proportional to the RMSE as $
    \lambda_k = \frac{\frac{1}{\text{RMSE}_k}}{\sum_{k=1}^K\frac{1}{\text{RMSE}_k} }$ \cite{nowotarski2014empirical}.
The resulting RMSE of the combined forecast using this method is 0.379.
Notably, the RMSE of the value-oriented forecast obtained with PH (B) is higher with 0.382.
Yet, as we observe in row RMSE in Table~\ref{tab:results_binary}, the value-oriented forecasts systematically improve upon the RMSE method. 
PH (B) achieves a cost saving of \$33164.

\subsection{Scalabilty}
We test the scalability of the proposed method, by running PFPH (CR) for the 2736-bus Summer Peak Polish system to which we added 21 400MW wind farms. We used the same data for load and wind forecasts as described in Section~\ref{ssec:data_description} above.
Table~\ref{tab:results_scalibility} summarizes the results. The training time remains manageable. Even for one full year of historical data, the algorithm converges after about 20h and improves the outcome in testing.

\begin{table}
\centering
\setlength{\tabcolsep}{3pt}
\caption{Results of PFPH for 2736 bus system} 
\vspace{-6pt}
\footnotesize
\renewcommand{\arraystretch}{1.2} 
\begin{tabular}{lccccc}
\hline
\text{Time Period}& $\lambda_1^*$&$\lambda_2^*$&\text{Time (s)}&\text{TST* [\$] }& $\Delta_c$ \\
\hline
one month &0.002&0.998& 17640& 117440&  -1167\\
one year &0.643&0.356& 73740 &113965& -606 \\
\hline
\end{tabular}
\label{tab:results_scalibility}
\end{table}

\section{conclusion}
We presented a method for value-oriented forecast combinations using progressive hedging (PH), unlocking high-fidelity, at-scale models and large-scale datasets in training.
We derived a one-shot reference model and discussed its scaling issues and presented the proposed PH approach alongside a modification that further reduces computation time.
Our case study demonstrated the effectiveness of value-oriented forecast combinations and showed the scalability of the the proposed method. 
Unit commitment and real-time dispatch cost were reduced by 1.8\% on average and we were able to obtain forecast combinations for the 2736 Polish system using a full year of historical data within 20 hours.
The method presented in this paper unlocks follow-up research on more context-aware forecast combination models as well as options to train models that provide advanced insights, such as forecast purchasing decisions.

{\color{black}
\appendix
\section{Convex relaxation of binary variables}
\label{appendix:a}
UC-R replaces the feasible set for each generator with its convex hull by relaxing $u_{g,t,d}$ to $[0,1]$ and by adding valid inequalities to tighten the relaxation of time-coupling constraints $\forall g \in [G]$ \cite{hua2016convex}:
\renewcommand{\theequation}{A.\arabic{equation}}
\setcounter{equation}{0}
\allowdisplaybreaks
\begin{subequations}
\begin{align}
& p_{g,t-1,d} \! \leq \overline{R}_g u_{g,t-1,d} \!+\!\! (\overline{P}_g\! -\! \overline{R}_g)(u_{g,t,d}\! -\! y_{g,t,d}) \ \ \forall t \!\in\! [2, T] \label{eq:ramprelax1}\\
&   p_{g,t,d} \! \leq \!\overline{P}_g u_{g,t,d} \!- \!(\overline{P}_g \!-\! \overline{R}_g\!) y_{g,t,d} \!\quad \forall t\!\! \in\! [2, T],\! \\
&  p_{g,t,d} - p_{g,t-1,d}  \leq (\underline{P}_g + R_g) u_{g,t,d} - \underline{P}_g u_{ g,t-1,d} - \nonumber \\ &(\underline{P}_g + R_g -\overline{R}_g) y_{t,d,g} \quad \forall t \in [2, T]\\
&   p_{g,t-1,d} - p_{g,t,d} \leq \overline{R}_g u_{ g,t-1,d} - (\overline{R}_g - R_g) u_{g,t,d} - \nonumber \\& (\underline{P}_g + R_g -\overline{R}_g) y_{g,t,d} \quad \forall t \in [2, T] \label{eq:ramprelax2}\\
&  u_{ g,t,d} \ge 0 \quad \forall t \in [2,T]\label{eq:endrelax}
\end{align}
\end{subequations}
\allowdisplaybreaks[0]
So, the resulting UC-R is:
\begin{equation}
\min \eqref{objectiveucrt} \text{ s.t. } \{
\eqref{eq:UC_balance} - \eqref{eq:UC_cur}, \eqref{eq:ramprelax1} - \eqref{eq:endrelax}
\} .\label{mainrelxappendix}
\end{equation}
}

\bibliographystyle{IEEEtran}
\bibliography{literature}

\begin{thebibliography}{10}
\providecommand{\url}[1]{#1}
\csname url@samestyle\endcsname
\providecommand{\newblock}{\relax}
\providecommand{\bibinfo}[2]{#2}
\providecommand{\BIBentrySTDinterwordspacing}{\spaceskip=0pt\relax}
\providecommand{\BIBentryALTinterwordstretchfactor}{4}
\providecommand{\BIBentryALTinterwordspacing}{\spaceskip=\fontdimen2\font plus
\BIBentryALTinterwordstretchfactor\fontdimen3\font minus
  \fontdimen4\font\relax}
\providecommand{\BIBforeignlanguage}[2]{{%
\expandafter\ifx\csname l@#1\endcsname\relax
\typeout{** WARNING: IEEEtran.bst: No hyphenation pattern has been}%
\typeout{** loaded for the language `#1'. Using the pattern for}%
\typeout{** the default language instead.}%
\else
\language=\csname l@#1\endcsname
\fi
#2}}
\providecommand{\BIBdecl}{\relax}
\BIBdecl

\bibitem{carriere2019integrated}
T.~Carriere \emph{et~al.}, ``An integrated approach for value-oriented energy
  forecasting and data-driven decision-making application to renewable energy
  trading,'' \emph{IEEE Trans. Smart Grid}, vol.~10, no.~6, 2019.

\bibitem{stratigakos2024decision}
A.~Stratigakos \emph{et~al.}, ``Decision-focused linear pooling for
  probabilistic forecast combination,'' \emph{International Journal of
  Forecasting}, 2024.

\bibitem{morales2023prescribing}
J.~M. Morales \emph{et~al.}, ``Prescribing net demand for two-stage electricity
  generation scheduling,'' \emph{Oper. Res. Perspect.}, vol.~10, 2023.

\bibitem{donti2017task}
P.~Donti \emph{et~al.}, ``Task-based end-to-end model learning in stochastic
  optimization,'' \emph{Adv. Neural Inf. Process. Syst.}, vol.~30, 2017.

\bibitem{zhang2024toward}
Y.~Zhang \emph{et~al.}, ``Toward value-oriented renewable energy forecasting:
  An iterative learning approach,'' \emph{IEEE Trans. Smart Grid}, 2024.

\bibitem{dias2025application}
J.~Dias~Garcia \emph{et~al.}, ``Application-driven learning: A closed-loop
  prediction and optimization approach applied to dynamic reserves and demand
  forecasting,'' \emph{Oper. Res.}, vol.~73, no.~1, 2025.

\bibitem{mieth2024prescribed}
R.~Mieth \emph{et~al.}, ``Prescribed robustness in optimal power flow,''
  \emph{Electric Power Systems Research}, vol. 235, 2024.

\bibitem{wang2023forecast}
X.~Wang \emph{et~al.}, ``Forecast combinations: An over 50-year review,''
  \emph{International Journal of Forecasting}, vol.~39, no.~4, 2023.

\bibitem{roald2023power}
L.~A. Roald \emph{et~al.}, ``Power systems optimization under uncertainty: A
  review of methods and applications,'' \emph{Electric Power Systems Research},
  vol. 214, p. 108725, 2023.

\bibitem{hua2016convex}
B.~Hua \emph{et~al.}, ``A convex primal formulation for convex hull pricing,''
  \emph{IEEE Trans. Power Syst.}, vol.~32, no.~5, 2016.

\bibitem{scholtes2001convergence}
S.~Scholtes, ``Convergence properties of a regularization scheme for
  mathematical programs with complementarity constraints,'' \emph{SIAM J.
  Optim}, vol.~11, no.~4, 2001.

\bibitem{ghazanfariharandi2025value}
M.~Ghazanfariharandi \emph{et~al.}, ``Value-oriented forecast combinations for
  unit commitment,'' \emph{arXiv preprint arXiv:2503.13677}, 2025.

\bibitem{kleinert2021survey}
T.~Kleinert \emph{et~al.}, ``A survey on mixed-integer programming techniques
  in bilevel optimization,'' \emph{EURO Journal on Computational Optimization},
  vol.~9, p. 100007, 2021.

\bibitem{rockafellar1991scenarios}
R.~T. Rockafellar \emph{et~al.}, ``Scenarios and policy aggregation in
  optimization under uncertainty,'' \emph{Math. Oper. Res.}, vol.~16, no.~1,
  1991.

\bibitem{gade2016obtaining}
D.~Gade \emph{et~al.}, ``Obtaining lower bounds from the progressive hedging
  algorithm for stochastic mixed-integer programs,'' \emph{Mathematical
  Programming}, vol. 157, pp. 47--67, 2016.

\bibitem{OpenPowerSystemDataplatform}
\BIBentryALTinterwordspacing
{Open Power System Data Platform}. [Online]. Available:
  \url{https://data.open-power-system-data.org/}
\BIBentrySTDinterwordspacing

\bibitem{dicopoulos2021weather}
J.~Dicopoulos \emph{et~al.}, ``Weather research and forecasting model
  validation with nrel specifications over the new york/new jersey bight for
  offshore wind development,'' in \emph{OCEANS 2021}.\hskip 1em plus 0.5em
  minus 0.4em\relax IEEE, 2021.

\bibitem{RUCOOL2019}
\BIBentryALTinterwordspacing
{RUCOOL}. (2019) Rutgers weather research and forecasting model. Accessed:
  2025-03-10. [Online]. Available:
  \url{https://tds.marine.rutgers.edu/thredds/dodsC/cool/ruwrf/wrf_4_1_3km_processed/WRF_4.1_3km_Processed_Dataset_Best.html}
\BIBentrySTDinterwordspacing

\bibitem{now23}
\BIBentryALTinterwordspacing
N.~Bodini \emph{et~al.} (2020) 2023 national offshore wind data set (now-23).
  [Online]. Available: \url{https://data.openei.org/submissions/4500}
\BIBentrySTDinterwordspacing

\bibitem{gaertner2020iea}
E.~Gaertner \emph{et~al.}, ``Iea wind tcp task 37: definition of the iea
  15-megawatt offshore reference wind turbine,'' National Renewable Energy
  Lab.(NREL), Golden, CO (United States), Tech. Rep., 2020.

\bibitem{lubin2023jump}
M.~Lubin \emph{et~al.}, ``Jump 1.0: Recent improvements to a modeling language
  for mathematical optimization,'' \emph{Math. Program. Comput.}, vol.~15,
  no.~3, 2023.

\bibitem{nowotarski2014empirical}
J.~Nowotarski \emph{et~al.}, ``An empirical comparison of alternative schemes
  for combining electricity spot price forecasts,'' \emph{Energy Economics},
  vol.~46, 2014.

\end{thebibliography}
\clearpage

\section*{Supplementary Material}

\subsection{Details on UC-R}
\label{appendix:a}
Typically, a UC problem determines a schedule of startup, shutdown, and power production for a set of power generation units. Constraints come in two forms: Coupling constraints enforce system-wide requirements (e.g., power balance); Individual (or private) constraints enforce technical limits of individual units. 
These constraints define the feasible set of each unit as well as the operational flexibility of the system.

In UC-R, the individual feasible set for each generating unit is replaced by its convex hull as proposed in \cite{hua2016convex}. 
In problem \eqref{eq:mainderetministicuc}, the private constraints that apply to each $g \in [G]$ are state-transition constraint \eqref{eq:UC_UC3}, minimum up/down time constraints \eqref{eq:UC_UC1}--\eqref{eq:UC_UC2}, dispatch level limits \eqref{eq:UC_powerlimit}, and ramping constraints \eqref{eq:UC_ramping1}--\eqref{eq:UC_ramping2}.

Without the time-coupling ramping constraints, the unit’s feasible set of the UC problem (with binary variables) is:
\[
\mathcal{X}_g \!=\! \left\{ \bm p_g \in \mathbb{R}^{T}, \bm y_{g}    \in \mathbb{R}^{T}, \bm u_{g}    \in \{0,1\}^{T-1}\ \!\!\middle|\ \! \!\text{\eqref{eq:UC_UC1}--\eqref{eq:UC_powerlimit}, \eqref{eq:uasbinary}} \right\}.
\]
The convex hull of this set $\mathcal{X}_g$ is:
\[
\text{conv}(\mathcal{X}_g) = \left\{ \bm p_g \in \mathbb{R}^{T}, \bm y_{g}    \in \mathbb{R}^{T}, \bm u_{g} \in \mathbb{R}^{T-1} \middle|\ \text{\eqref{eq:UC_UC1}--\eqref{eq:UC_powerlimit}} \right\}.
\]
We now consider ramping constraints \eqref{eq:UC_ramping1}--\eqref{eq:UC_ramping2}. When ramping constraints are included in the definition of $\mathcal{X}_g$,
the
valid inequalities \eqref{eq:ramprelax1}--\eqref{eq:ramprelax2} as shown above can be applied to any two consecutive time periods to tighten
the approximation of $\text{conv}(\mathcal{X}_g)$ \cite{hua2016convex}.
\allowdisplaybreaks
So  for each unit $g \in [G]$ in our  formulation, the
feasible region is defined to be:
\[
 \left\{ \bm p_g \!\in \mathbb{R}^{T}, \bm y_{g}    \in \mathbb{R}^{T}, \bm u_{g} \in \mathbb{R}^{T-1} \middle|\ \!\! \text{\eqref{eq:UC_UC1}--\eqref{eq:UC_ramping2}, \eqref{eq:ramprelax1}--\eqref{eq:endrelax}} \! \right\}.
\]
which is a tractable approximation of $\text{conv}(\mathcal{X}_g)$. So the resulting UC-R is defined as
\renewcommand{\theequation}{S.\arabic{equation}}
\setcounter{equation}{0}
\begin{subequations}
\label{mainrelxx}
\begin{align}
 \min \quad 
&  \eqref{objectiveucrt}\nonumber\\
\text{s.t.} \quad
    & \eqref{eq:UC_balance} - \eqref{eq:UC_flow2}, \quad \eqref{eq:UC_shed} - \eqref{eq:UC_cur}\\
    & \left\{ \bm p_g , \bm y_{g}    , \bm u_{g}  \right\} \in \text{conv}(\mathcal{X}_g) \quad \forall g \in [G].
\end{align}
\end{subequations}
(We note that \eqref{mainrelxx} and \eqref{mainrelxappendix} are identical).

\subsection{KKT conditions}
This section shows the resulting Karush-Kuhn-Tucker (KKT) conditions for the lower-level unit commitment relaxation (UC-R) problem in \eqref{eq:mainbilevel} for the exemplary case that minimum up- and downtimes are $\underline{\ell}_g=\overline{\ell}_g=1$. 
\allowdisplaybreaks
\begin{subequations}
\begin{align}
& \forall d \in [D]:\nonumber\\
    & (\hat{f}_{l,t,d}): = - \sum_{i} A_{l,i} \nu^1_{i,t,d} - \nu^2_{l,t,d} - \Psi^1_{l,t,d} + \Psi^2_{l,t,d} =0 \nonumber \\&\qquad \forall l \in [L] , \forall t \in [T] \label{eq:stationarity1}\\
   &(\hat{\theta}_{i,t,d}): \sum_{l} A_{l,i} B_l \nu^2_{l,t,d} = 0 \quad \forall i \in [N/N_0] , \forall t \in [T]\\
   & (\hat{\theta}_{i,t,d})\!:\!\! \sum_{l} A_{l,i} B_l \nu^2_{l,t,d} - \nu^3_{ref,t,d} = 0 \quad \forall i \in [N_0] , \forall t \in [T]\\
   & \forall g \in [G]:\nonumber\\
   &
   \left\{ \begin{array}{l}
(u_{g,t,d}): - \Psi^{6}_{g,t,d} \Bar{P}_g + \Psi^{7}_{g,t,d} \underline{P}_g + c^{SD}_g - \Psi^{5}_{g,t+1,d} -\nonumber\\\\\Psi^{8}_{g,t+1,d}  R_g + \Psi^{8}_{g,t+1,d} \Bar{R}_g -\! \Psi^{10}_{g,t+1,d} \Bar{R}_g \!+\! \Psi^{12}_{g,t+1,d} \underline{P}_g \nonumber\\\\- \Psi^{13}_{g,t+1,d} \Bar{R}_g + \Psi^{4}_{t+\underline{\ell}_g,g,d} - \Psi^{14}_{g,t,d} = 0 \quad t= 1 \\\\
  (u_{g,t,d}): -c^{SD}_g + \Psi^{5}_{g,t,d} -  \Psi^{6}_{g,t,d} \Bar{P}_g +  \Psi^{7}_{g,t,d} \underline{P}_g  - \nonumber\\\\ \Psi^{9}_{g,t,d} (R_g + \Bar{R}_g) - \Psi^{10}_{g,t,d} (\Bar{P}_g - \Bar{R}_g) -  \Psi^{12}_{g,t,d}  \nonumber\\\\ (\underline{P}_g + R_g)+\Psi^{13}_{g,t,d} (\Bar{R}_g + R_g) + c^{SD}_g - \Psi^{5}_{g,t+1,d} -  \nonumber\\\\\Psi^{8}_{g,t+1,d} R_g +\! \Psi^{8}_{g,t+1,d} \Bar{R}_g -\!\Psi^{10}_{g,t+1,d} \Bar{R}_g +\! \Psi^{12}_{g,t+1,d} \underline{P}_g \nonumber\\\\ - \Psi^{13}_{g,t+1,d} \Bar{R}_g -\Psi^{3}_{g,t,d} + \Psi^{4}_{t+\underline{\ell}_g,g,d} - \Psi^{14}_{g,t,d}= 0 \quad \nonumber\\\\ \overline{\ell}_g +1 \le t \le T-1\\\\
  (u_{g,t,d}): -c^{SD}_g + \Psi^{5}_{g,t,d} -  \Psi^{6}_{g,t,d} \Bar{P}_g +  \Psi^{7}_{g,t,d} \underline{P}_g  - \nonumber\\\\ \Psi^{9}_{g,t,d} (R_g + \Bar{R}_g) - \Psi^{10}_{g,t,d} (\Bar{P}_g - \Bar{R}_g) -  \Psi^{12}_{g,t,d}  \nonumber\\\\(\underline{P}_g + R_g)+ \Psi^{13}_{g,t,d}   - \Psi^{3}_{g,t,d} - \Psi^{14}_{g,t,d}= 0 \quad t= T
\end{array} \right.\\
&
\left\{ \begin{array}{l}
 (p_{g,t,d}): C^g - \nu^1_{i(g),t,d}  + \Psi^{6}_{g,t,d} -  \Psi^{7}_{g,t,d} -  \Psi^{8}_{g,t,d} + \nonumber\\\\ \Psi^{9}_{g,t+1,d} +  \Psi^{10}_{g,t+1,d} -\!  \Psi^{12}_{g,t+1,d} +\!  \Psi^{13}_{g,t+1,d} =\! 0 \quad t= 1 \\\\
 (p_{g,t,d}): C^g - \nu^1_{i(g),t,d}  + \Psi^{6}_{g,t,d} -  \Psi^{7}_{g,t,d} +  \Psi^{8}_{g,t,d} - \nonumber\\\\\Psi^{9}_{g,t,d} +  \Psi^{11}_{g,t,d} +  \Psi^{12}_{g,t,d} -  \Psi^{13}_{g,t,d} -  \Psi^{8}_{g,t+1,d} +  \nonumber\\\\\Psi^{9}_{g,t+1,d} +  \Psi^{10}_{g,t+1,d}- \Psi^{12}_{g,t+1,d} + \Psi^{13}_{g,t+1,d}= 0\nonumber\\\\ \quad 2\le t \le T-1 \\\\
  (p_{g,t,d}): C^g - \nu^1_{i(g),t,d}  + \Psi^{6}_{g,t,d} -  \Psi^{7}_{g,t,d} +  \Psi^{8}_{g,t,d}  \nonumber\\\\-  \Psi^{9}_{g,t,d}+  \Psi^{11}_{g,t,d} + \Psi^{12}_{g,t,d} -  \Psi^{13}_{g,t,d} = 0 \quad t= T
\end{array} \right.\\
&
\left\{ \begin{array}{l}
 (y_{g,t,d}): c^{SU}_g + c^{SD}_g - \Psi^{5}_{g,t,d} + \Psi^{10}_{g,t,d} \Bar{P}_g - \Psi^{10}_{g,t,d} \Bar{R}_g \nonumber\\\\ + \Psi^{11}_{g,t,d} \Bar{P}_g - \Psi^{11}_{g,t,d} \Bar{R}_g + \Psi^{12}_{g,t,d} (\underline{P}_g + R_g - \Bar{R}_g) + \nonumber\\\\\Psi^{13}_{g,t,d} (\underline{P}_g + R_g - \Bar{R}_g) -\Psi^{15}_{g,t,d}= 0 \quad t= 1 \\\\
 (y_{g,t,d}):c^{SU}_g + c^{SD}_g - \Psi^{5}_{g,t,d} + \Psi^{10}_{g,t,d} \Bar{P}_g - \Psi^{10}_{g,t,d} \Bar{R}_g \nonumber\\\\ +\Psi^{11}_{g,t,d} \Bar{P}_g - \Psi^{11}_{g,t,d} \Bar{R}_g + \Psi^{12}_{g,t,d} (\underline{P}_g + R_g - \Bar{R}_g) +\nonumber\\\\ \Psi^{13}_{g,t,d} (\underline{P}_g + R_g - \Bar{R}_g) + \Psi^{3}_{g,t,d} + \Psi^{4}_{g,t,d} - \Psi^{15}_{g,t,d}=\! 0 \nonumber\\\\\quad  2 \le t \le T
\end{array} \right.\\
   &(\hat{l}_{i,t,d}): c_i^{shed} - \nu_{i,t,d} - \psi^{17} + \psi^{18}=0 \quad\forall i \in [N] ,\nonumber \\ &\quad\forall t \in [T] \label{eq:stationarity2}\\
   &(\hat{w}_{i,t,d}): c_i^{cur} + \nu_{i,t,d} - \psi^{19} + \psi^{20}=0 \quad\forall i \in [N] ,\nonumber \\ &\quad\forall t \in [T] \label{eq:stationarity3}\\
   & \eqref{eq:UC_balance} - \eqref{eq:UC_cur}, \quad \eqref{eq:ramprelax1} - \eqref{eq:endrelax}\quad \text{[Primal feasibility]} \label{eq:primal feasibility} \\
&  \Psi_{g,t,d}^j f_j(p_{g,t,d}, u_{g,t,d}, y_{g,t,d} ) =0  \quad  \forall g \in [G], \forall j \in [J],\nonumber \\ & \forall t \in [T]\label{eq:cnonlinearslackness1}\\
    & \Psi_{i,t,d}^b f_b(\hat{l}_{i,t,d}, \hat{w}_{i,t,d} ) =0  \quad\forall i \in [N], \forall b \in [B], \forall t \in [T] \label{eq:cnonlinearslackness2} \\
    &  \Psi_{l,t,d}^q f_q(\hat{f}_{l,t,d} ) =0 \quad   \forall l \in [L], \forall q \in [Q], \forall t \in [T] \label{eq:cnonlinearslackness3}\\
&  \Psi^b_{i,t,d} \ge 0 \quad  \forall i \in [N], \forall b \in [B], \forall t \in [T]  \label{eq:dualfeasibility1}\\
& \Psi^j_{g,t,d} \ge 0 \quad  \forall g \in [G], \forall j \in [J], \forall t \in [T] \\
&\Psi^q_{l,t,d} \ge 0 \quad  \forall l \in [L], \forall q \in [Q], \forall t \in [T] \\
& \nu_{i,t,d}:free \label{eq:dualfeasibility}
\end{align}
\end{subequations}
where $\bm{\Psi}^b$, $\bm{\Psi}^j$, and $\bm{\Psi}^q$ are the Lagrange multipliers of the  inequality constraints on each node $i$ represented by $f_b(\cdot)$, each generator $g \in [G]$ represented by $f_j(\cdot)$, and each transmission line $l$ represented by $f_q(\cdot)$.
Values $B$, $Q$, and $J$ denote the respective numbers of constraints/dual variables.
$A_{l,i}$ is the incidence matrix of transmission lines: 1 or –1 indicates that node $i$ is the from/to node for line $l$, and 0 otherwise.
Also, $\bm{\nu}$ is the Lagrange multiplier related to equality constraints. The constraints \eqref{eq:stationarity1}-\eqref{eq:stationarity3}, \eqref{eq:primal feasibility}, \eqref{eq:cnonlinearslackness1}-\eqref{eq:cnonlinearslackness3}, and \eqref{eq:dualfeasibility1}-\eqref{eq:dualfeasibility} are, respectively, the stationarity, primal feasibility, complementary slackness, and dual feasibility conditions of the lower-level problem in \eqref{eq:mainbilevel}. 
These KKT conditions are necessary for formulating the single-level problem. We then replace these KKT conditions with UC-R in \eqref{eq:mainbilevel}.

\subsection{Resolving KKT nonlinearity}

This section shows details on how methods for resolving the nonlinearities of the KKT conditions have been applied.

First, we can address the resulting non-linearity in the complementarity slackness conditions \eqref{eq:cnonlinearslackness1}--\eqref{eq:cnonlinearslackness3} using a regularization approach from \cite{scholtes2001convergence} which yields a local optimum \cite{scholtes2001convergence}. This method replaces \eqref{eq:cnonlinearslackness1}-\eqref{eq:cnonlinearslackness3} with:
\allowdisplaybreaks
\begin{subequations}
\begin{align}
& \forall d \in [D], \forall t \in [T] : \nonumber \\
    & \sum_{j \in [J]} \Psi_{g,t,d}^j f_j(p_{g,t,d}, u_{g,t,d}, y_{g,t,d} ) \le \epsilon  \quad \forall g \in [G] \label{eq:nonlinearslackness1}\\
    &\sum_{b \in [B]} \Psi_{i,t,d}^b f_b(\hat{l}_{i,t,d}, \hat{w}_{i,t,d} ) \le \epsilon  \quad \forall i \in [N] \label{eq:nonlinearslackness2} \\
    & \sum_{q \in [Q]} \Psi_{l,t,d}^q f_q(\hat{f}_{l,t,d} ) \le \epsilon \quad  \forall l \in [L]. \label{eq:nonlinearslackness3}
\end{align}%
\end{subequations}%
\allowdisplaybreaks[0]
Here, $\epsilon$ represents a small non-negative scalar that enables the reformulation of the KKT condition into a parametrized nonlinear problem that can be solved by modern off-the-shelf non-linear solvers. 
As outlined in the main body of the paper, we denote this method as ST-N. 

Alternatively, the complementarity slackness conditions can be linearized using Fortuny–Amat (``Big-M''):
\allowdisplaybreaks
\begin{subequations}
\begin{align}
& \forall t \in [T], \forall d \in [D] : \nonumber\\
    & 0\le \Psi^j_{g,t,d} \le M z^j_{g,t,d} \quad  \forall g \in [G],\ \forall j \in [J] \label{eq:linearslackness1}\\
& 0\le f_j(\hat{f}_{g,t,d} ) \le M (1-z^j_{g,t,d}) \quad \forall g \in [G],\ \forall j \in [J] \\
    & 0\le \Psi^b_{i,t,d} \le M z^b_{i,t,d} \quad  \forall i \in [N],\ \forall b \in [B] \\
& 0\le f_b(\hat{f}_{i,t,d} ) \le M (1-z^b_{i,t,d}) \quad \forall i \in [N],\ \forall b \in [B]\\
    & 0\le \Psi^q_{l,t,d} \le M z^q_{l,t,d} \quad \forall l \in [L],\ \forall q \in [Q] \\
& 0\le f_q(\hat{f}_{l,t,d} ) \le M (1-z^q_{l,t,d}) \quad \forall l \in [L],\ \forall q \in [Q] \label{eq:linearslackness2}
\end{align}%
\end{subequations}%
\allowdisplaybreaks[0]%
where $\bm{z}$ are binary variables, and $M \in \mathbb{R}^+$ is a large enough constant. As outlined in the main body of the paper we denote this method as ST-M. 

The resulting single-level equivalent of \eqref{eq:mainbilevel} is:
\begin{align}
 \min_{\lambda_1,...,\lambda_K} \quad 
&  \eqref{ew:bilevel_obj}\nonumber\\
\text{s.t.} \quad
    & \eqref{eq: forecast combination} - \eqref{eq:secondstageinbilevel}, \quad \eqref{eq:stationarity1} - \eqref{eq:primal feasibility}, \quad \eqref{eq:dualfeasibility1} - \eqref{eq:dualfeasibility} \label{eq:monoliticproblem}\\
    & \eqref{eq:nonlinearslackness1}-\eqref{eq:nonlinearslackness3}\text{ for ST-N} \,\text{ \textbf{or} }\nonumber\\& \, \eqref{eq:linearslackness1}-\eqref{eq:linearslackness2}\text{ for ST-M}.\nonumber 
\end{align}
The training problem in \eqref{eq:monoliticproblem} solves a two-level network-constrained problem for each time $t$ over all days $d$ and includes binary (ST-M) or non-linear (ST-N) structures.

\vspace{12pt}

\end{document}